\documentclass[reqno,12pt]{amsart}

\usepackage{amsmath, amsthm, amssymb, amsfonts, amsthm}
\usepackage{amsxtra, amscd, geometry, graphicx,epic,eepic}
\usepackage{mathrsfs}
\usepackage[all]{xy}

\newtheorem{thm}{Theorem}[section]
\newtheorem{lem}[thm]{Lemma}
\newtheorem{prop}[thm]{Proposition}
\newtheorem{cor}[thm]{Corollary}

\theoremstyle{definition}

\theoremstyle{remark}
\newtheorem{remark}[thm]{Remark}

\theoremstyle{plain}

\newcommand{\Z}{{\mathbb{Z}}}

\newcommand{\Q}{\mathbb{Q}}
\newcommand{\R}{{\mathbb{R}}}

\newcommand{\N}{{\mathbb{N}}}

\newcommand{\FF}{{\mathcal{F}}}
\newcommand{\TT}{\mathscr T }
\newcommand{\SSS}{\mathscr S}

\newcommand{\DD}{{\mathscr{D}}}
\newcommand{\RR}{\mathscr{R}}
\newcommand{\eps}{\varepsilon}

\newcommand{\area}{{\operatorname{Area}}}

\newcommand{\ii}{\mathrm{i}}

\newcommand{\EE}{\mathscr E}
\newcommand{\NN}{\mathscr N}
\newcommand{\MM}{\mathscr M}
\newcommand{\WW}{\mathscr W}

\newcommand{\per}{\operatorname{per}}
\newcommand{\eper}{\operatorname{eper}}
\newcommand{\tr}{\operatorname{tr}}
\newcommand{\Tr}{\operatorname{Tr}}

\numberwithin{equation}{section}

\begin{document}

\title[]{Products of matrices
$\left[ \begin{matrix} 1 & 1 \\ 0 & 1
\end{matrix}\right]$ and
$\left[ \begin{matrix} 1 & 0 \\ 1 & 1 \end{matrix} \right]$ and
the distribution of reduced quadratic irrationals}

\author[]{Florin P. Boca}

\address{Department of Mathematics, University of Illinois at Urbana-Champaign,
1409 W. Green Street, Urbana IL 61801, USA}
\address{and}
\address{Institute of Mathematics "Simion Stoilow" of the Romanian
Academy, P.O. Box 1-764, RO-014700 Bucharest, Romania}

\address{E-mail: fboca@math.uiuc.edu}

\date{}

\begin{abstract}
Let $\Phi(N)$ denote the number of products of matrices $\left[
\begin{smallmatrix} 1 & 1 \\ 0 & 1 \end{smallmatrix}\right]$ and
$\left[ \begin{smallmatrix} 1 & 0 \\ 1 & 1
\end{smallmatrix}\right]$ of trace equal to $N$, and
$\Psi(N)=\sum_{n=3}^N \Phi(n)$ be the number of such products of
trace between $3$ and $N$.
We prove an asymptotic formula of type $\Psi(N)=c_1 N^2 \log
N+c_2N^2 +O_\eps(N^{7/4+\eps})$ as $N\rightarrow \infty$. As a
result, the Dirichlet series $\sum_{n=3}^\infty \Phi(n) n^{-s}$
has a meromorphic extension in the half-plane $\Re
(s)>\frac{7}{4}$ with a single, order two pole at $s=2$. Our
estimate also improves on an asymptotic result of Faivre
concerning the distribution of reduced quadratic irrationals,
providing an explicit upper bound for the error term.
\end{abstract}


\maketitle

\section{Introduction}
Consider the matrices
\begin{equation*}
A=\left[ \begin{matrix} 1 & 0 \\ 1 & 1 \end{matrix} \right] \qquad
\mbox{\rm and} \qquad B=\left[ \begin{matrix} 1 & 1 \\ 0 & 1
\end{matrix} \right]
\end{equation*}
and the (free) multiplicative monoid $\MM$ they generate.
The problem of estimating
\begin{equation*}
\Phi(N) =\# \{ C\in \MM \, :\, \operatorname{Tr}(C)= N\}
\end{equation*}
and
\begin{equation*} \Psi(N) =\sum\limits_{n=3}^N \Phi(n)=\# \{
C\in \MM \, :\, 3 \leq \operatorname{Tr}(C)\leq N\}
\end{equation*}
for large $N$ came across in the study of a number-theoretic spin
chain model in statistical mechanics introduced in \cite{KO}, and
further investigated in \cite{FKO} and \cite{CKK}. In \cite{KOPS}
the estimate
\begin{equation}\label{1.1}
\Psi (N)=\frac{N^2 \log N}{\zeta(2)}+O(N^2 \log \log N) \qquad
\qquad  (N\rightarrow \infty)
\end{equation}
was proved, using $L$-functions and a result from \cite{Fa}
concerning the distribution of reduced quadratic irrationals.

In this paper we improve the estimate \eqref{1.1}. Our approach
relies on a result concerning the distribution of multiplicative
inverses, which is a consequence of Weil's bound on Kloosterman
sums.

\medskip

\begin{thm}\label{T1.1} We have
\begin{equation}\label{1.2}
\Psi(N)=c_1 N^2 \log N +c_2 N^2 +\Psi_0(N),
\end{equation}
where
\begin{equation}\label{1.3}
\Psi_0 (N)\ll_\eps N^{\frac{7}{4}+\eps} \quad\qquad (N\rightarrow
\infty),
\end{equation}
and
\begin{equation*}
c_1 =\frac{1}{\zeta(2)} , \qquad c_2 =\frac{1}{\zeta(2)} \bigg(
\gamma-\frac{3}{2} -\frac{\zeta^\prime(2)}{\zeta(2)}\bigg).
\end{equation*}
\end{thm}

Using the Mellin transform representation of Dirichlet series and
changing $N$ to $N-2$ in the right-hand side of \eqref{1.2} we
obtain for $\Re (s)>2$
\begin{equation*}
\begin{split}
Z(s) & =\sum\limits_{n=3}^\infty \Phi (n)n^{-s}=s\int_3^\infty
\Psi(x) (x-2)^{-s-1} \, dx
\\ & =s\int_1^\infty (c_1 x^2\log x+c_2 x^2)x^{-s-1}\,
dx+s\int_1^\infty \Psi_0 (x)x^{-s-1}\, dx \\ &
=\frac{2c_1}{(s-2)^2}+\frac{c_1+2c_2}{s-2}+c_2+s\int\limits_1^\infty
\Psi_0(x) x^{-s-1}\, dx.
\end{split}
\end{equation*}
Since the function $\Psi_0(x)$ satisfies the growth condition
\eqref{1.3} the integral above converges for $\Re
(s)>\frac{7}{4}$, and thus the right-hand side defines an analytic
continuation of $Z(s)$ to the half-plane $\Re(s)>\frac{7}{4}$ with
the point $s=2$ removed.

The contribution to the main term in the asymptotic formula above
only comes from words of odd length in $\mathscr{M}$. We estimate
the contribution $\Psi_{\operatorname{ev}}(N)$ to $\Psi(N)$ of words of
even length which begin in $A$ and end in $B$, proving

\medskip

\begin{prop}\label{P1.2}
$\qquad \displaystyle \Psi_{\operatorname{ev}}(N) =\frac{N^2\log
2}{2\zeta(2)}+O_\eps (N^{\frac{7}{4}+\eps})\qquad \qquad
(N\rightarrow \infty).$
\end{prop}

This extends \cite[Prop.4.5]{KOPS} which states that
\begin{equation*}
\Psi_{\operatorname{ev}}(N) \sim \frac{N^2\log 2}{2\zeta(2)} ,
\end{equation*}
without any estimate on the error term.

Using a transfer operator associated with the Gauss map, Fredholm
theory, and Ikehara's tauberian theorem, Faivre \cite{Fa} proved
the asymptotic formula
\begin{equation*}
\sum\limits_{\rho (\omega)<X} 1 \sim \frac{\operatorname{e}^X \log
2}{2\zeta(2)}\qquad \qquad (X\rightarrow \infty)
\end{equation*}
for the number of reduced quadratic irrationals $\omega$ of length
$\rho(\omega)$ at most $X$. Since the final argument relies on a
tauberian theorem, no explicit bound was found for the error term.
In the last section we use Proposition \ref{P1.2} and an explicit
identification from \cite{KOPS} between products of matrices $A$
and $B$ starting with $B$ and ending with $A$, and reduced
quadratic irrationals, to prove
\begin{prop}\label{P1.3}
$\quad \displaystyle \sum\limits_{\rho (\omega)<X} 1 =
\frac{\operatorname{e}^X\log 2}{2\zeta(2)} +O_\eps
(\operatorname{e}^{(\frac{7}{8}+\eps)X})\qquad \qquad (X\rightarrow
\infty).$
\end{prop}

\section{Products of $A$'s and $B$'s and continued fractions}
If
\begin{equation*}
\alpha=[a_1,a_2,\dots]=\frac{1}{a_1+\frac{1}{a_2+\frac{1}{a_3+\frac{1}{\dots}}}}
\end{equation*}
is a reduced continued fraction with positive integers $a_i$, the
$k^{\operatorname{th}}$ convergent
\begin{equation*}
\frac{p_k}{q_k}=[a_1,\dots,a_k]
\end{equation*}
is given by pairs $(p_n,q_n)$ of relatively prime integers defined
recursively as
\begin{equation}\label{2.1}
\begin{cases}
p_0=0,\ p_1=1,\ p_n=a_n p_{n-1}+p_{n-2} ,\\
q_0=1,\ q_1=a_1,\ q_n=a_n q_{n-1}+q_{n-2} ,
\end{cases}
\end{equation}
and satisfying $0\leq p_n\leq q_n$ and the equality
\begin{equation*}
p_{n-1} q_n-p_n q_{n-1}=(-1)^n.
\end{equation*}

If $(p_n)$ and $(q_n)$ satisfy \eqref{2.1} for every $0\leq n\leq
k$, then
\begin{equation*}
\frac{q_{n-1}}{q_n}=\frac{q_{n-1}}{a_n q_{n-1}+q_{n-2}}
=\frac{1}{a_n+\frac{q_{n-2}}{q_{n-1}}}  ,
\end{equation*}
showing that
\begin{equation}\label{2.2}
\frac{q_{k-1}}{q_{k}} =[a_{k},\dots,a_1].
\end{equation}

Consider the matrices
\begin{equation*}
J=\left[ \begin{matrix} 0 & 1 \\ 1 & 0 \end{matrix} \right] \qquad
\mbox{\rm and} \qquad M(a)=\left[ \begin{matrix} a & 1 \\ 1 & 0
\end{matrix} \right] .
\end{equation*}

If $(p_n,q_n)$ is as in \eqref{2.1} then, as noticed in
\cite{KOPS}, we have
\begin{equation}\label{2.3}
M(a_1)\dots M(a_n)=\left[ \begin{matrix} q_n & q_{n-1} \\ p_n &
p_{n-1} \end{matrix} \right].
\end{equation}
When combined with
\begin{equation*}
B^k A^\ell =M(k)M(\ell),\qquad k,\ell \in \Z,
\end{equation*}
equality \eqref{2.3} yields
\begin{equation}\label{2.4}
B^{a_1} A^{a_2} \cdots B^{a_{2m-1}} A^{a_{2m}} =M(a_1)\cdots
M(a_{2m})=\left[ \begin{matrix} q_{2m} & q_{2m-1} \\ p_{2m} &
p_{2m-1} \end{matrix} \right] .
\end{equation}
From \eqref{2.1}, \eqref{2.4}, and
\begin{equation*}
B=A^T=JAJ,\qquad A=JBJ,
\end{equation*}
we also infer that
\begin{equation}\label{2.5}
\begin{split}
A^{a_1} B^{a_2}\cdots A^{a_{2m-1}}B^{a_{2m}} & =J\left[
\begin{matrix} q_{2m} & q_{2m-1} \\ p_{2m} & p_{2m-1} \end{matrix}
\right] J=\left[ \begin{matrix} p_{2m-1} & p_{2m} \\ q_{2m-1} &
q_{2m} \end{matrix} \right], \\
B^{a_1}A^{a_2} \cdots A^{a_{2m}} B^{a_{2m+1}} & = \left[
\begin{matrix} q_{2m} & q_{2m-1} \\ p_{2m} & p_{2m-1} \end{matrix} \right]
\left[ \begin{matrix} 1 & a_{2m+1} \\ 0 & 1 \end{matrix} \right] =
\left[ \begin{matrix} q_{2m} & q_{2m+1} \\ p_{2m} & p_{2m+1}
\end{matrix} \right] , \\
A^{a_1} B^{a_2} \cdots B^{a_{2m}} A^{a_{2m+1}} & =\left[
\begin{matrix} p_{2m+1} & p_{2m} \\ q_{2m+1} & q_{2m}
\end{matrix}\right] .
\end{split}
\end{equation}
All matrices in the products from \eqref{2.4} and \eqref{2.5} have
determinant $1$. We denote
\begin{equation*}
\begin{split}
\WW_{\operatorname{ev}}(N) & =\{ (a_1,\dots,a_{2m})\in \N^{2m} \,
:\,m\geq 1,\
\Tr(B^{a_1}A^{a_2} \cdots B^{a_{2m-1}}A^{a_{2m}})\leq N\},\\
\WW_{\operatorname{odd}}(N) & =\{ (a_1,\dots,a_{2m+1})\in \N^{2m+1} \,
:\,m\geq 1,\  \Tr(B^{a_1}A^{a_2} \cdots
A^{a_{2m}}B^{a_{2m+1}})\leq N\}.
\end{split}
\end{equation*}

We consider the sets
\begin{equation*}
\SSS_{\operatorname{ev}}(N)=\left\{ \left[ \begin{matrix}
q^\prime & q
\\ p^\prime & p \end{matrix} \right] \, :\, \begin{matrix} 0\leq
p\leq q,\ 0\leq p^\prime \leq q^\prime ,\ q^\prime >q,\\
p+q^\prime \leq N,\ pq^\prime -p^\prime q=1 \end{matrix}\right\}
\end{equation*}
and
\begin{equation*}
\SSS_{\operatorname{odd}}(N)=\left\{ \left[ \begin{matrix} q & q^\prime \\
p & p^\prime
\end{matrix} \right] \, :\,
\begin{matrix} 0\leq p\leq q,\ 0\leq p^\prime \leq q^\prime,\ q^\prime \geq q,\\
p^\prime+q\leq N,\ p^\prime q-pq^\prime =1 \end{matrix} \right\},
\end{equation*}
of cardinality $\Psi_{\operatorname{ev}}(N)$ and respectively
$\Psi_{\operatorname{odd}}(N)$, and the maps defined as
\begin{equation*}
\begin{split}
\beta_{\operatorname{ev}}(a_1,\dots,a_{2m}) & =
B^{a_1} A^{a_2} \cdots B^{a_{2m-1}} A^{a_{2m}}=M(a_1)\cdots M(a_{2m}), \\
\beta_{\operatorname{odd}}(a_1,\dots,a_{2m+1}) & =
B^{a_1} A^{a_2} \cdots A^{a_{2m}} B^{a_{2m+1}} = M(a_1)\cdots M(a_{2m+1}) J,
\end{split}
\end{equation*}
from $\bigcup_1^\infty \WW_{\operatorname{ev}}(N)$ to $\bigcup_1^\infty
\SSS_{\operatorname{ev}}(N)$, and respectively from $\bigcup_1^\infty
\WW_{\operatorname{odd}}(N)$ to $\bigcup_1^\infty \SSS_{\operatorname{odd}}(N)$.
As a consequence of \eqref{2.2}, $\beta_{\operatorname{ev}}$ and
$\beta_{\operatorname{odd}}$ are injective. It follows from \eqref{2.4}
and \eqref{2.5} that $\beta_{\operatorname{ev}}( \WW_{\operatorname{ev}}(N)
)\subseteq \SSS_{\operatorname{ev}}(N)$ and $\beta_{\operatorname{odd}}(
\WW_{\operatorname{odd}}(N)) \subseteq \SSS_{\operatorname{odd}}(N)$. To check
$\SSS_{\operatorname{ev}}(N) \subseteq
\beta_{\operatorname{ev}}(\WW_{\operatorname{ev}}(N))$,
let $\left[ \begin{smallmatrix} q^\prime & q \\
p^\prime & p \end{smallmatrix} \right] \in \SSS_{\operatorname{ev}}(N)$.
With $K=\big[ \frac{q^\prime}{q}\big]$ we have $0\leq q^\prime
-Kq<q$, $0\leq p^\prime -Kp$, and $(p^\prime -Kp)q-(q^\prime
-Kq)p=1$. Since $q^\prime >q$ are relatively prime, we also have
$K\leq \frac{q^\prime -1}{q}\leq \frac{q^\prime -p^\prime}{q-p}$,
and thus $p^\prime -Kp\leq q^\prime -Kq$. Since
\begin{equation*}
\left[ \begin{matrix} q^\prime & q \\ p^\prime & p \end{matrix}
\right] \left[ \begin{matrix} 0 & 1 \\ 1 & -K \end{matrix} \right]
=\left[ \begin{matrix} q^\prime & q \\ p^\prime & p \end{matrix}
\right] M(K)^{-1}=
\left[ \begin{matrix} q & q^\prime -Kq \\
p & p^\prime -Kp
\end{matrix} \right],
\end{equation*}
it follows, by replacing $(q,q^\prime)$ by $(q^\prime -Kq,q)$ and
$(p,p^\prime)$ by $(p^\prime -Kp,p)$ and performing this process
until $q^\prime -Kq$ becomes equal to $1$, that
the matrix $\left[ \begin{smallmatrix} q^\prime & q \\
p^\prime & p \end{smallmatrix} \right]$ is written as a product of
$k$ matrices of form $M(K)$. Since $q^\prime >q$, $k$ ought to be
even and therefore $\left[
\begin{smallmatrix} q^\prime & q \\ p^\prime & p \end{smallmatrix} \right]
\in \beta_{\operatorname{ev}} ( \WW_{\operatorname{ev}}(N))$. One shows in a
similar way that $\SSS_{\operatorname{odd}}(N) \subseteq
\beta_{\operatorname{odd}}( \WW_{\operatorname{odd}}(N))$.

This proves that the elements of $\mathscr{M}$ are uniquely
represented as products of $A$'s and $B$'s. It also implies that
\begin{equation}\label{2.6}
\Psi(N)=2\Psi_{\operatorname{ev}} (N)+2\Psi_{\operatorname{odd}}(N).
\end{equation}

\section{Estimating $\Psi_{\operatorname{ev}} (N)$}
To estimate $\Psi_{\operatorname{ev}}(N)$, we first keep $q^\prime$ and
$q$ constant. From $pq^\prime -p^\prime q=1$ and $q<q^\prime$, it
follows that $q$ and $q^\prime$ are relatively prime, and that $p$
is uniquely determined as $p=\overline{q^\prime}$, where
$\overline{q^\prime}$ is the unique integer in $\{ 1,\dots,q\}$
for which $q^\prime \overline{q^\prime}=1 \pmod{q}$. It is obvious
that $p^\prime :=\frac{pq^\prime -1}{q}\leq q^\prime$ and the map
\begin{equation*}
\left\{ (q,q^\prime)\, : \, \begin{matrix} q<q^\prime \leq N,\
(q,q^\prime)=1\\  q^\prime +\overline{q^\prime} \leq N
\end{matrix} \right\} \ni (q,q^\prime) \ \mapsto \
\left[ \begin{matrix} q^\prime & q \\
p^\prime =\frac{\overline{q^\prime} q^\prime -1}{q} &
p=\overline{q^\prime}
\end{matrix} \right]\in \mathscr{S}_{\operatorname{ev}}(N)
\end{equation*}
is a one-to-one correspondence. Replacing $q^\prime$ by $y$ and
$\overline{q^\prime}$ by $x$ we can write
\begin{equation}\label{3.1}
\Psi_{\operatorname{ev}} (N)= \sum\limits_{\substack{q<q^\prime \leq N \\
\overline{q^\prime}+q^\prime \leq N}} 1=\sum\limits_{q<N}
\sum\limits_{\substack{q<y\leq N \\ 0<x\leq \min \{ q,N-y\} \\
xy=1 \hspace{-6pt} \pmod{q}}} 1.
\end{equation}

For each $y\in (0,N]$, there is at most one $x\in (0,q)$ such that
$xy=1\pmod{q}$; whence the trivial estimate
\begin{equation*}
\Psi_{\operatorname{ev}} (N)\ll N^2 .
\end{equation*}
To give a more precise estimate for $\Psi_{\operatorname{ev}}(N)$, we
shall define for $q>1$ integer and $\Omega$ subset in $\R^2$ the
number $\NN_q (\Omega)$ of (relatively prime) integers $(x,y)\in
\Omega$ such that $xy=1\pmod{q}$. It is known (see for instance
\cite{BCZ}) that Weil's bound on Kloosterman sums yields for any
intervals $I$ and $J$ of length at most $q$ the estimate
\begin{equation*}
\NN_q(I\times J)=\frac{\varphi(q)}{q^2}\ \vert I\vert \, \vert
J\vert +O_\varepsilon (q^{\frac{1}{2}+\varepsilon}).
\end{equation*}
This immediately extends to intervals of arbitrary size as
\begin{equation}\label{3.2}
\NN_q(I\times J)=\frac{\varphi(q)}{q^2}\ \vert I\vert \, \vert
J\vert +O_\varepsilon \bigg( q^{\frac{1}{2}+\eps} \Big(
1+\frac{\vert I\vert}{q}\Big)\Big( 1+\frac{\vert J\vert}{q} \Big)
\bigg).
\end{equation}
Another easy but useful consequence of \eqref{3.2} is given next.

\begin{lem}\label{L3.1}
Suppose that
\begin{equation*}
\Omega =\{ (x,y)\, :\, \alpha \leq x\leq \beta,\ f_1(x)\leq y\leq
f_2(x)\} ,
\end{equation*}
with $C^1$-functions $f_1 \leq f_2$ on $[\alpha,\beta]$. For every
positive integer $T$ we have
\begin{equation*}
\NN_q(\Omega)=\frac{\varphi(q)}{q^2} \, \area(\Omega)+\EE_q,
\end{equation*}
with
\begin{equation*}
\EE_q \ll_\eps \frac{\beta-\alpha}{Tq}\ \big(V_\alpha^\beta
(f_1)+V_\alpha^\beta(f_2)\big)+Tq^{\frac{1}{2}+\eps}\bigg(
1+\frac{\beta-\alpha}{Tq}\bigg) \bigg( 1+\frac{\| f_1\|_\infty+\|
f_2\|_\infty}{q}\bigg),
\end{equation*}
where $V_\alpha^\beta (f_i)$ denotes the total variation of $f_i$
on $[\alpha,\beta]$.
\end{lem}

\begin{proof} One can take without loss of generality $f_1=0$ and
$f_2 \geq 0$. Partitioning $[\alpha,\beta]$ into $T$ intervals
$I_i$ of equal size and denoting by $M_i$ and $m_i$ the maximum,
respectively the minimum, of $f_2$ on $I_i$, we clearly have
\begin{equation*}
\sum\limits_{i=0}^{T-1} \NN_q \big( I_i \times [ 0,m_i]\big) \leq
\NN_q (\Omega)\leq \sum\limits_{i=0}^{T-1} \NN_q \big( I_i \times
[ 0,M_i]\big).
\end{equation*}
By \eqref{3.2} we infer that
\begin{equation*}
\NN_q \big( I_i\times [0,m_i]\big) =\frac{\varphi(q)}{q^2} \cdot
\frac{\beta-\alpha}{T} \ m_i +O_\eps \bigg( q^{\frac{1}{2}+\eps}
\Big( 1+\frac{\beta-\alpha}{Tq}\Big)\Big( 1+\frac{\|
f_2\|_\infty}{q}\Big) \bigg).
\end{equation*}
The statement follows from this, the similar estimate for $\NN_q
( I_i\times [0,M_i])$, and from
\begin{equation*}
\begin{split}
\int_\alpha^\beta f_2(x)\, dx & =\sum\limits_{i=0}^{T-1}
\frac{\beta-\alpha}{T}\ m_i+O\bigg( \frac{\beta-\alpha}{T} \
V_\alpha^\beta (f_2)\bigg) \\ & =\sum\limits_{i=0}^{T-1}
\frac{\beta-\alpha}{T}\ M_i+O\bigg( \frac{\beta-\alpha}{T} \
V_\alpha^\beta (f_2)\bigg).
\end{split}
\end{equation*}
\end{proof}

Equality \eqref{3.1} can also be written as
\begin{equation}\label{3.3}
\Psi_{\operatorname{ev}} (N) =\sum\limits_{q<N} \NN_q (\Omega_{N,q}),
\end{equation}
where $\Omega_{N,q}=\{ (x,y)\, :\, q<y\leq N,\ 0<x\leq \min \{
q,N-y\} \}$ coincides with  the trapezoid $\Omega_{N,q}^{(1)}=\{
(x,y)\ :\, 0<x\leq q<y\leq N-x\}$ if $q\leq \frac{N}{2}$, and with
the triangle $\Omega_{N,q}^{(2)}=\{(x,y)\, :\, 0<x<N-q,\ q<y\leq
N-x\}$ if $q>\frac{N}{2}$. Employing \eqref{3.3} and applying
Lemma \ref{L3.1} to $\Omega_{N,q}^{(1)}$ and $\Omega_{N,q}^{(2)}$
we infer that
\begin{equation*}
\begin{split}
\Psi_{\operatorname{ev}} (N)& =\sum\limits_{q\leq \frac{N}{2}} \NN_q
(\Omega_{N,q}^{(1)})+\sum\limits_{\frac{N}{2}<q<N} \NN_q
(\Omega_{N,q}^{(2)}) \\
& =\sum\limits_{q\leq \frac{N}{2}} \frac{\varphi(q)}{q^2} \cdot
\frac{q(2N-3q)}{2}+\sum\limits_{\frac{N}{2}<q<N}
\frac{\varphi(q)}{q^2}\cdot \frac{(N-q)^2}{2}+\EE_1 (N)+\EE_2(N),
\end{split}
\end{equation*}
with
\begin{equation*}
\begin{split}
\EE_1(N) & \ll_\eps \sum\limits_{q\leq \frac{N}{2}} \bigg(
\frac{q}{Tq}\cdot q+Tq^{\frac{1}{2}+\eps} \cdot \frac{N}{q} \bigg)
\ll_\eps \frac{N^2}{T}+TN^{\frac{3}{2}+\eps},\\
\EE_2(N) & \ll_\eps \sum\limits_{\frac{N}{2}<q<N} \bigg(
\frac{N}{Tq}\cdot N+Tq^{\frac{1}{2}+\eps} \Big(
1+\frac{N}{Tq}\Big) \frac{N}{q}\bigg)\ll_\eps
\frac{N^2}{T}+TN^{\frac{3}{2}+\eps} .
\end{split}
\end{equation*}

Taking $T=[N^{\frac{1}{4}}]$ and using standard summation results
(see formulas in Corollary 4.5 and its proof below) this gives
\begin{equation*}
\Psi_{\operatorname{ev}}(N)=\frac{N^2 \log 2}{2\zeta(2)} +O_\eps (N^{\frac{7}{4}+\eps}),
\end{equation*}
which ends the proof of Proposition \ref{P1.2}.

\section{Estimating $\Psi_{\operatorname{odd}}(N)$}
To estimate $\Psi_{\operatorname{odd}} (N)$, we first keep $q$ and $p$
fixed. The general solution of $p^\prime q-pq^\prime =1$ is given
by
\begin{equation*}
p^\prime =\overline{q}+pt,\qquad q^\prime =\frac{p^\prime
q-1}{p}=\frac{q\overline{q}-1}{p}+qt,\qquad t\in \Z,
\end{equation*}
where $\overline{q}$ is the unique integer in $\{ 1,\dots,p\}$ for
which $\overline{q} q =1 \pmod{p}$. Since $p^\prime
>p$, one has $t\geq 1$. The map
\begin{equation*}
\left\{ (q,p,t)\, :\, \begin{matrix} 0\leq p<q,\ (p,q)=1\\ 1\leq
t\leq \big[ \frac{N-q-\overline{q}}{p} \big] \end{matrix} \right\}
\ni (q,p,t)\ \mapsto \ \left[ \begin{matrix} q & q^\prime
=\frac{q\overline{q}-1}{p}+qt \\ p & p^\prime =\overline{q}+pt
\end{matrix}\right] \in \SSS_{\operatorname{odd}} (N)
\end{equation*}
is a bijection. Replacing $p$ by $a$, $q$ by $y$, and
$\overline{q}$ by $x$, we can write
\begin{equation*}
\Psi_{\operatorname{odd}}(N)  =\sum\limits_{q<N}
\sum\limits_{\substack{p<q \\ q+\overline{q} \leq N}} \bigg[
\frac{N-q-\overline{q}}{p}\bigg]  =\sum\limits_{a<N}
\sum\limits_{\substack{a<y<N \\ 0<x\leq \min \{ a,N-y \} \\
xy=1\hspace{-6pt} \pmod{a}}} \bigg[ \frac{N-y-x}{a}\bigg].
\end{equation*}
When $N-y<a$ we have $N-y-x<a$, and the contribution of such terms
to $\Psi_{\operatorname{odd}}(N)$ is null. So we only consider $N-y\geq
a$. This gives $a<\frac{N}{2}$ and thus
\begin{equation*}
\Psi_{\operatorname{odd}} (N)=\sum\limits_{a<\frac{N}{2}}
\sum\limits_{\substack{a<y\leq N-a \\ 0<x\leq a \\
xy=1 \hspace{-6pt} \pmod{a}}} \left[ \frac{N-y-x}{a}\right] .
\end{equation*}

From
\begin{equation*}
\left[ \frac{N-y-x}{a}\right] \leq \frac{N-y-x}{a}< \frac{N-a}{a}
=\frac{N}{a}-1
\end{equation*}
it follows that
\begin{equation*}
\left[ \frac{N-y-x}{a} \right] \leq \left[ \frac{N}{a}\right] -1.
\end{equation*}
The set of points $(x,y)\in (0,a]\times (a,N-a)$ for which $\big[
\frac{N-x-y}{a}\big]=i$ coincides with
\begin{equation*}
\Omega_{N,a,i}:=\{ (x,y)\in (0,a]\times (a,N-a)\, :\, N-(i+1)a<x+y
\leq N-ia\}.
\end{equation*}
We can thus write
\begin{equation}\label{4.1}
\begin{split}
\Psi_{\operatorname{odd}}(N) & =\sum\limits_{a<\frac{N}{2}}
\sum\limits_{i=1}^{\big[ \frac{N}{a}\big]-1} i\NN_a
(\Omega_{N,a,i})
\\ & =\sum\limits_{a<\frac{N}{2}} \sum\limits_{j=1}^{\big[
\frac{N}{a}\big]-1} \big( j\NN_a (\Omega_{N,a,j}^{(1)}
)+(j-1)\NN_a (\Omega_{N,a,j}^{(2)})\big) ,
\end{split}
\end{equation}
where the sets
\begin{equation*}
\begin{split}
& \Omega_{N,a,j}^{(1)} = \big\{ (x,y)\in (0,a]\times \big(
N-(j+1)a,N-ja\big] ,\ x+y \leq N-ja\big\},\\ &
\Omega_{N,a,j}^{(2)} =\big\{ (x,y)\in (0,a]\times \big(
N-(j+1)a,N-ja\big],\
x+y>N-ja\big\} ,\quad 1\leq j\leq \bigg[ \frac{N}{a}\bigg]-2,\\
& \Omega^{(1)}_{N,a,[N/a]-1} =\left\{ (x,y)\, :\,
\begin{matrix}
0<x\leq N-\big[ \frac{N}{a}\big] a,\ a<y\leq N-\big( \big[
\frac{N}{a}\big] -1\big)a,\\ x+y\leq N-\big( \big[
\frac{N}{a}\big] -1\big) a \end{matrix}\right\},
\\ & \Omega^{(2)}_{N,a,[N/a]-1} =\left\{ (x,y)\, :\,
\begin{matrix}
0<x\leq a,\ a<y\leq N-\big( \big[ \frac{N}{a}\big] -1\big)a,\\
x+y> N-\big( \big[ \frac{N}{a}\big] -1\big) a
\end{matrix}\right\},
\end{split}
\end{equation*}
give a partition of the trapezoid $\{ (x,y)\, :\, 0<x\leq a <y\leq
N-a-x\}$. Since $\NN_a (\Omega)$ does not change when $\Omega$ is
translated by integer multiples of $(a,0)$, the right-hand side in
\eqref{4.1} can also be expressed as
\begin{equation*}
\sum\limits_{j=1}^{\big[ \frac{N}{a}\big]-2} \NN_a (\SSS_{N,a,j})\
+\NN_a (\SSS_{N,a,[N/a]-1}),
\end{equation*}
where for $1\leq j\leq \big[ \frac{N}{a}\big]-2$ we set
\begin{equation*}
\SSS_{N,a,j} =\big\{ (x,y)\in (0,ja] \times \big( N-(j+1)a,
N-ja\big],\ x+y\leq N-a\big\},\end{equation*} and for $j=\big[
\frac{N}{a}\big]-1$ we set
\begin{equation*} \SSS_{N,a,[N/a]-1} =\left\{ (x,y)\, :\, 0<x\leq
N-2a,\ a<y\leq N-\bigg( \bigg[ \frac{N}{a}\bigg]-1\bigg) a,\
x+y\leq N-a\right\} .
\end{equation*}
The sets $\SSS_{N,a,j}$, $1\leq j\leq \big[ \frac{N}{a}\big]-1$,
are mutually disjoint and their union is the triangle
\begin{equation*}
\TT_{N,a}=\{ (x,y)\, :\, 0<x\leq N-2a,\ a<y\leq N-a-x\},
\end{equation*}
and thus we get
\begin{equation*}
\Psi_{\operatorname{odd}}(N)=\sum\limits_{a<\frac{N}{2}} \NN_a
(\TT_{N,a}).
\end{equation*}
Since $\NN_a(\Omega)$ is invariant under translations by $(0,a)$,
this further gives
\begin{equation}\label{4.2}
\Psi_{\operatorname{odd}}(N)=\sum\limits_{a<\frac{N}{2}} \NN_a
(\widetilde{\TT}_{N,a}),
\end{equation}
where
\begin{equation*}
\widetilde{\TT}_{N,a} =\{ (x,y)\, :\, 0<x\leq N-2a,\ 0<y\leq
N-2a-x\}
\end{equation*}
is the translated triangle $\TT_{N,a} -(0,a)$.

\begin{lem}\label{L4.1}
For every $0<c<1$ we have
\begin{equation*}
\sum\limits_{a\leq N^c} \NN_a (\widetilde{\TT}_{N,a})
=\sum\limits_{a\leq N^c} \frac{\varphi(a)}{a^2} \cdot
\frac{(N-2a)^2}{2}+O(N^{1+c}).
\end{equation*}
\end{lem}

\begin{proof} Set $K=\big[ \frac{N}{a}\big]-2\geq 1$. We partition
the triangle $\widetilde{\TT}_{N,a}$ as $\DD_{N,a}\cup \RR_{N,a}$,
with
\begin{equation*}
\DD_{N,a} =\bigcup\limits_{i=1}^K \ (0,ia] \times \big(
(K-i)a,(K-i+1)a\big]
\end{equation*}
and $\RR_{N,a}=\widetilde{\TT}_{N,a} \setminus \DD_{N,a}$. As
$\DD_{N,a}$ is the union of $1+2+\cdots +K=\frac{K(K+1)}{2}$
disjoint squares of size $a$, we have
\begin{equation}\label{4.3}
\NN_a (\DD_{N,a})=\frac{K(K+1)}{2}\
\varphi(a)=\frac{\varphi(a)}{a^2}\ \area \DD_{N,a} .
\end{equation}

On the other hand it is clear that
\begin{equation*}
\NN_a (\RR_{N,a})\leq (K+1)\varphi(a) \leq \frac{N}{a}\cdot a=N,
\end{equation*}
which gives in turn
\begin{equation}\label{4.4}
\sum\limits_{a\leq N^c} \NN_a (\RR_{N,a}) \leq N^{1+c} .
\end{equation}
We also have
\begin{equation}\label{4.5}
\sum\limits_{a\leq N^c} \frac{\varphi(a)}{a^2}\ \area \RR_{N,a}
\leq \sum\limits_{a\leq N^c} \frac{1}{a}\ \area \RR_{N,a} \leq
\sum\limits_{a\leq N^c} \frac{1}{a}\cdot (K+1)\cdot \frac{a^2}{2}
\leq \sum\limits_{a\leq N^c} \frac{N}{a}\cdot a\leq N^{1+c}.
\end{equation}
Employing \eqref{4.4}, \eqref{4.3} and \eqref{4.5}, we gather
\begin{equation*}
\begin{split}
\sum\limits_{a\leq N^c} \NN_a (\widetilde{\TT}_{N,a}) &
=\sum\limits_{a\leq N^c} \NN_a (\DD_{N,a}) +O(N^{1+c})
=\sum\limits_{a\leq N^c} \frac{\varphi(a)}{a^2} \cdot \area
\DD_{N,a} +O(N^{1+c}) \\
& =\sum\limits_{a\leq N^c} \frac{\varphi(a)}{a^2} \, (\area
\widetilde{\TT}_{N,a} -\area \RR_{N,a})+O(N^{1+c}) \\
& =\sum\limits_{a\leq N^c} \frac{\varphi(a)}{a^2} \cdot \area
\widetilde{\TT}_{N,a} +O(N^{1+c}) \\ & =\sum\limits_{a\leq N^c}
\frac{\varphi(a)}{a^2} \cdot \frac{(N-2a)^2}{2} +O(N^{1+c}). \qedhere
\end{split}
\end{equation*} 
\end{proof}

\begin{lem}\label{L4.2}
For every $0<c<1$ and every integer $T>1$ we have
\begin{equation*}
\sum\limits_{N^c <a<\frac{N}{2}} \NN_a (\widetilde{\TT}_{N,a})
=\sum\limits_{N^c <a<\frac{N}{2}} \frac{\varphi(a)}{a^2} \cdot
\frac{(N-2a)^2}{2}+O_\eps \bigg( \frac{N^2 \log
N}{T}+TN^{\frac{3}{2}+\eps} +N^{2+(\eps-\frac{1}{2})c}\bigg) .
\end{equation*}
\end{lem}

\begin{proof}   Applying Lemma \ref{L3.1} we get
\begin{equation*}
\begin{split}
\NN_a (\widetilde{\TT}_{N,a}) & =\frac{\varphi(a)}{a^2}\cdot \area
(\widetilde{\TT}_{N,a}) +O_\eps \bigg( \frac{N}{Ta}\cdot
N+Ta^{\frac{1}{2}+\eps} \Big( 1+\frac{N}{Ta}\Big) \cdot
\frac{N}{a}\bigg) \\
& =\frac{\varphi(a)}{a^2} \cdot \frac{(N-2a)^2}{2} +O_\eps \bigg(
\frac{N^2}{Ta}+TNa^{-\frac{1}{2}+\eps}+N^2
a^{-\frac{3}{2}+\eps}\bigg),
\end{split}
\end{equation*}
which is summed in the range $a\in (N^c,\frac{N}{2})$ to get
\begin{equation*}
\begin{split}
\sum\limits_{N^c <a<\frac{N}{2}} & \frac{\varphi(a)}{a^2}\cdot
\frac{(N-2a)^2}{2}+O_\eps \bigg( \frac{N^2 \log
N}{T}+TN\sum\limits_{a=1}^N a^{-\frac{1}{2}+\eps}+N^2
\sum\limits_{a>N^c} a^{-\frac{3}{2}+\eps}\bigg) \\
& = \sum\limits_{N^c <a<\frac{N}{2}} \frac{\varphi(a)}{a^2}\cdot
\frac{(N-2a)^2}{2}+O_\eps \bigg( \frac{N^2 \log
N}{T}+TN^{\frac{3}{2}+\eps} +N^{2+(\eps-\frac{1}{2})c}\bigg). \qedhere
\end{split}
\end{equation*}
\end{proof}

Taking $T=[ N^{\frac{1}{4}}]$ and (any) $c\in
[\frac{1}{2},\frac{3}{4}]$, the previous two lemmas, together with
\eqref{4.2}, yield

\begin{cor}\label{C4.3}
$\quad \displaystyle
\Psi_{\operatorname{odd}}(N)=\sum\limits_{a<\frac{N}{2}}
\frac{\varphi(a)}{a^2}\cdot \frac{(N-2a)^2}{2}+O_\eps
(N^{\frac{7}{4}+\eps}).$
\end{cor}

\begin{lem}\label{L4.4}
$\quad S_N:=\displaystyle \sum\limits_{a<\frac{N}{2}}
\frac{\varphi(a)(N-2a)^2}{2a^2} =C_N+O(N)$, where
\begin{equation*}
C_N=\frac{N^2}{2\zeta(2)} \left( \log N+\gamma -\log
2-\frac{3}{2}-\frac{\zeta^\prime(2)}{\zeta(2)}\right).
\end{equation*}
\end{lem}

\begin{proof}
Employing the Dirichlet series
\begin{equation*}
\sum\limits_{a=1}^\infty \frac{\varphi(a)}{a^s}
=\frac{\zeta(s-1)}{\zeta(s)}  ,\qquad s=\sigma+it,\ \sigma >2,
\end{equation*}
and the Perron integral formula $(\sigma_0>0$)
\begin{equation*}
\frac{1}{2\pi \ii} \int_{\sigma_0-\ii \infty}^{\sigma_0+\ii
\infty} \frac{y^s}{s(s+1)(s+2)} \, ds =\begin{cases} 0 & \mbox{\rm if
$0\leq y\leq 1$,} \\ \frac{1}{2}\big(1-\frac{1}{y}\big)^2 &
\mbox{\rm if $y\geq 1$,}
\end{cases}
\end{equation*}
with $y=\frac{N}{2a}$ we infer that
\begin{equation*}
S_N=\frac{1}{2\pi\ii}
\int_{\sigma_0-\ii\infty}^{\sigma_0+\ii\infty} g(s)\ ds,
\end{equation*}
with
\begin{equation*}
\begin{split}
g(s) & =\frac{N^{s+2}}{2^s (s+1)(s+2)\zeta(s+2)}\cdot
\frac{\zeta(s+1)}{s} \\ & =\frac{N^{s+2}}{2^s(s+1)(s+2)\zeta(s+2)}
\left(\frac{1}{s^2}+\frac{\gamma}{s}+O(1)\right)\qquad
(s\rightarrow 0).
\end{split}
\end{equation*}
In the region $\Re s>-2$ the function $g$ is meromorphic with a
removable singularity at $s=-1$ and a pole $C_N=h^\prime (0)$ at
$s=0$, where
\begin{equation*}
h(s)=\frac{N^{s+2}(1+\gamma s)}{2^s(s+1)(s+2)\zeta(s+2)}.
\end{equation*}
A direct calculation gives
\begin{equation*}
C_N=h^\prime (0)=\frac{N^2}{2\zeta(2)} \left( \log N+\gamma -\log
2-\frac{3}{2}-\frac{\zeta^\prime(2)}{\zeta(2)}\right).
\end{equation*}

We seek to change the contour of integration from
$\sigma=\sigma_0$ to the contour $\Gamma$ consisting of the five
line segments $s=\sigma_0\pm it$ $(t\geq T)$, $s=\sigma\pm \ii T$
$(-1\leq \sigma\leq \sigma_0)$, $s=-1+\ii t$ $(-T\leq t\leq T)$,
getting
\begin{equation}\label{4.6}
\frac{1}{2\pi\ii} \int_{\sigma_0-\ii\infty}^{\sigma_0+\ii\infty}
g(s)\, ds=C_N+\frac{1}{2\pi\ii} \int_\Gamma g(s)\, ds.
\end{equation}
It remains to show that the contribution of the integral on
$\Gamma$ is small. Note first that $\vert \zeta
(\overline{s})\vert=\vert \zeta (s)\vert$ gives $\vert
g(\overline{s})\vert=\vert g(s)\vert$. As a result only the case
$\Im s\geq 0$ will be considered next. Using standard estimates on
$\zeta$ (cf., e.g., \cite{Edw},\cite{Ing},\cite{Iv}) we have
\begin{equation}\label{4.7}
\int_T^\infty \vert g(\sigma_0+\ii t)\vert\,  dt\ll_{\sigma_0}
\int_T^\infty \frac{N^{2+\sigma_0}\log t}{2^{\sigma_0}t^3}\, dt
\ll_{\sigma_0,\eps} \frac{N^{2+\sigma_0}}{T^{2-\eps}}
\end{equation}
and
\begin{equation}\label{4.8}
\begin{split}
\int_{-1}^{\sigma_0} \vert g(\sigma+\ii T)\vert, d\sigma & \ll
\int_{-1}^{\sigma_0} \frac{N^{2+\sigma_0}}{T^3}\cdot \vert \zeta
(1+\sigma+\ii T)\vert\cdot\frac{1}{\vert\zeta (2+\sigma+\ii
T)\vert} \, d\sigma
\\ & \ll_{\sigma_0}
\frac{N^{2+\sigma_0}}{T^3} (T^{1/2}\log T ) \log^7
T\ll_{\sigma_0,\eps}\frac{N^{2+\sigma_0}}{T^{5/2-\eps}}.
\end{split}
\end{equation}
To estimate the contribution of the integrand on the segment
$-1+\ii t$ $(0<t\leq T)$ we follow closely the argument from
\cite{Hall}, pp. 216-217. The functional equation
\begin{equation*}
\frac{\zeta(\ii t)}{\zeta(1+\ii t)} =\chi (\ii t)\cdot \frac{\zeta
(1-\ii t)}{\zeta (1+\ii t)}\, ,\qquad \chi
(s)=\frac{(2\pi)^s}{2\Gamma (s)\cos \frac{\pi s}{2}}\, ,
\end{equation*}
and the equality
\begin{equation*}
\vert \Gamma (\ii t)\vert^2 =\frac{\pi}{t\sinh \pi t}
\end{equation*}
yield
\begin{equation*}
\left| \frac{\zeta (\ii t)}{\zeta (1+\ii t)}\right| =\vert \chi
(\ii t)\vert = \frac{1}{2\sqrt{\frac{\pi}{t\sinh \pi t}} \cdot
\cosh \frac{\pi t}{2}} =\sqrt{\frac{t\tanh \frac{\pi
t}{2}}{2\pi}}.
\end{equation*}
Employing also $\tanh t\leq t$ $(t\geq 0)$, we get
(independently of $T\geq \frac{2}{\pi}$)
\begin{equation}\label{4.9}
\int_0^T \vert g(-1+\ii t)\vert \, d t 
=\frac{N}{2} \int_0^T \bigg| \frac{\zeta (it)}{\zeta (1+it)}\bigg| \frac{dt}{t(1+t^2)} 
\ll \frac{N}{2} \int_0^T \frac{t}{t(1+t^2)}\, dt\ll N .
\end{equation}
The estimates \eqref{4.6}-\eqref{4.9} with, say, $T=N^2$ conclude
the proof.
\end{proof}

Note also

\begin{cor}\label{C4.5}
$\quad \displaystyle \sum\limits_{a<N}
\frac{\varphi(a)}{a^2}=\frac{1}{\zeta(2)} \bigg( \log
N+\gamma-\frac{\zeta^\prime(2)}{\zeta(2)}\bigg) +O \bigg(
\frac{\log N}{N}\bigg).$
\end{cor}

\begin{proof}
This is a consequence of Lemma \ref{L4.4} and of the well known
formulas
\begin{equation*} \sum\limits_{a<N}
\varphi(a)=\frac{N^2}{2\zeta(2)}+O(N\log N), \qquad
\sum\limits_{a<N} \frac{\varphi(a)}{a}=\frac{N}{\zeta(2)}+O(\log
N).
\end{equation*}
\end{proof}

Theorem \ref{T1.1} now follows from \eqref{2.6}, Proposition
\ref{P1.2} and Lemma \ref{L4.4}.

\begin{figure}[ht]
\includegraphics*[scale=.9, bb=80 10 330 140]{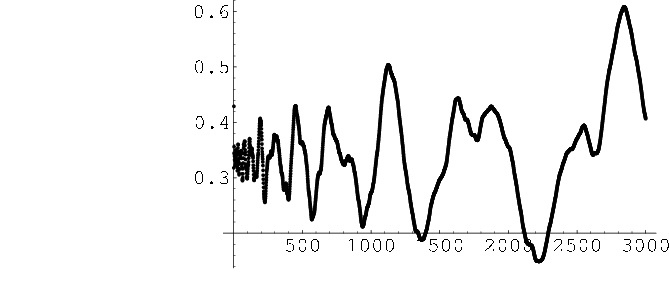}
\includegraphics*[scale=.9, bb=90 10 330 140]{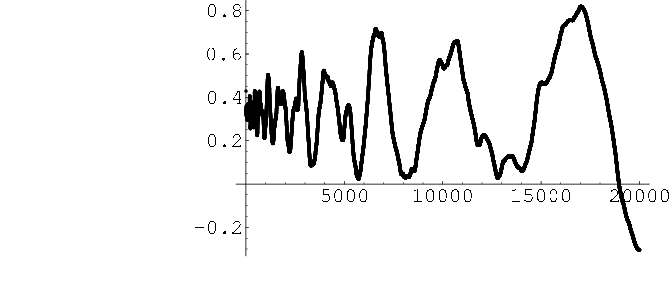}
\caption{\sl The plot of $S_N-C_N$, $N\leq 3000$, and respectively
$N\leq 20000$.} \label{Figure1}
\end{figure}

\begin{remark}\label{R4.6}
Numerical computations show that the error $O(N)$ given by Lemma
\ref{L4.4} on
\begin{equation*}
S_N-C_N =\sum_{a<\frac{N}{2}}
\frac{\varphi(a)(N-2a)^2}{2a^2}-\frac{N^2}{2\zeta(2)} \left( \log
N+\gamma-\frac{3}{2}-\log
2-\frac{\zeta^\prime(2)}{\zeta(2)}\right)
\end{equation*}
may not be optimal. Moreover, the graph of $S_N-C_N$ exhibits a
surprising regularity (cf. Figure \ref{Figure1}). This was brought
to our attention by the referee, who also kindly provided the
Mathematica notebook. One could hope to improve the theoretical
estimate of the error by shifting the segment $s=-1+\ii t$ further
left to the line $\Re s=-1-\delta$. The problem however is that
the argument $2+s$ will enter the critical strip $0<\Re s<1$ where
lower bound estimates for $\zeta$ are problematic.
\end{remark}

\begin{figure}[ht]
\includegraphics*[scale=0.95, bb=100 10 350 140]{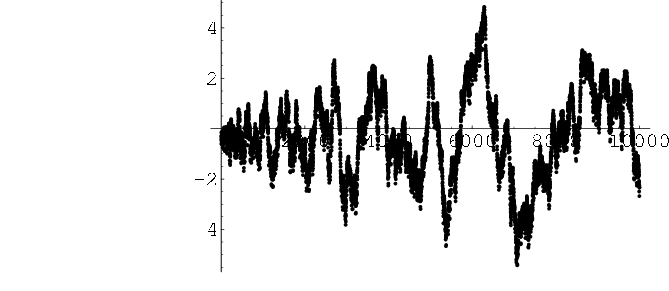}
\caption{\sl The plot of $\sum_{a<N/2}
\frac{\varphi(a)(N-2a)}{a}-\frac{N^2}{4\zeta(2)}$, $N\leq 10000$.}
\label{Figure2}
\end{figure}

\section{An application to the distribution of quadratic irrationals}

Let $x$ be an irrational number in $(0,1)$ with continued fraction
$[a_0(x),a_1(x),\dots ]$ and $\frac{p_n(x)}{q_n(x)}$ be its
$n^\mathrm{th}$ convergent. A classical result of P. L\' evy
states for almost every $x$ that
\begin{equation*}
\lim\limits_{n\rightarrow \infty} \frac{1}{n} \log
q_n(x)=\frac{\pi^2}{12\log 2}  .
\end{equation*}

The Gauss map $T:[0,1)\rightarrow [0,1)$ defined as $T(0)=0$ and
$T(x)=\frac{1}{x}-\big[ \frac{1}{x}\big]$ if $x\neq 0$, has the
well-known properties
\begin{equation*}
\begin{split}
T\big( [a_1,a_2,\dots ]\big) & =[a_2,a_3,\dots ],\\
 xT(x)\cdots T^{n-1}(x) & =(-1)^n \big( xq_n(x)-p_n(x)\big).
\end{split}
\end{equation*}

When $\omega$ is a quadratic irrational, it is well-known that the
limit $\beta(\omega)$ of $\frac{1}{n}\log q_n(\omega)$ exists, and
is called the L\' evy constant of $\omega$. Let $AX^2+BX+C$ be the
minimal integer polynomial of $\omega$ and $\Delta=B^2-4AC$. The
length of $\omega$ is defined as $\rho (\omega)=2\log
\eps_0(\omega)$, where $\eps_0(\omega)=\frac{1}{2}(u_0+v_0
\sqrt{\Delta})$ is the fundamental solution of the Pell equation
$u^2-\Delta v^2=4$.

We are interested in the set $\RR$ of all purely periodic
quadratic irrationals, aiming to evaluate
\begin{equation*}
\pi_0 (X)=\sum\limits_{\substack{\omega \in \RR \\ \rho
(\omega)<X}} 1 \qquad \qquad (X\rightarrow \infty).
\end{equation*}

Following \cite{KOPS}, one defines for each such
$\omega=\overline{[a_1,\dots,a_n]}$ with $n=\per (\omega)$ the
quantities
\begin{equation*}
\begin{split}
\eper(\omega) & =\begin{cases} n, & \mbox{\rm if $n=\per(\omega)$
even},\\ 2n, & \mbox{\rm if $n=\per(\omega)$ odd,} \end{cases}
\\ M(\omega) & =M(a_1)\cdots M(a_n),\\
\widetilde{M} (\omega) & =\Omega^+ =\begin{cases} M(\omega), &
\mbox{\rm if $n$ even,} \\ M(\omega)^2 , & \mbox{\rm if $n$ odd.}
\end{cases}
\end{split}
\end{equation*}
According to \cite{Fa}, Proposition 2.2, we have
\begin{equation*}
\begin{split}
& \eps_0(\omega)=\omega T(\omega) T^2(\omega)\cdots
T^{\eper(\omega)-1} (\omega),\\ & \rho (\omega)=2\log
\eps_0(\omega).
\end{split}
\end{equation*}

Writing $\Delta=f^2 \Delta_0$ for some fundamental discriminant
$\Delta_0$ and some positive integer $f$, one considers the group
of units of $\mathscr{O}_\Delta$
\begin{equation*}
E_\Delta=\bigg\{ \frac{u+v\sqrt{\Delta}}{2} \, :\, (u,v)\in \Z^2,\
u^2-\Delta v^2 =\pm 4 \bigg\}
\end{equation*}
with fundamental unit $\eps_\Delta>1$, in the quadratic field
$\mathfrak{K}=\Q (\sqrt{\Delta_0})$ endowed with $\Q$-valued norm
$\NN$ and trace $\tr$. One also considers the subgroup
\begin{equation*}
E_\Delta^+=\bigg\{ \frac{u+v\sqrt{\Delta}}{2} \, :\, (u,v)\in
\Z^2,\ u^2-\Delta v^2 =4 \bigg\}
\end{equation*}
of the totally positive units, which is generated by
\begin{equation*}
\eps^+_\Delta=\begin{cases} \eps_\Delta & \mbox{\rm if
$\NN(\eps_\Delta)=+1,$}\\
\eps^2_\Delta & \mbox{\rm if $\NN(\eps_\Delta)=-1.$} \end{cases}
\end{equation*}

In \cite{KOPS}, Section 2, the explicit isomorphism
\begin{equation*}
\lambda_\omega :\FF_\omega \rightarrow E_\Delta,\qquad
\lambda_\omega \left( \left[ \begin{matrix} a & b \\ c & d
\end{matrix} \right]\right) =c\omega+d
\end{equation*}
between the fixed point group
\begin{equation*}
\FF_\omega =\{ g\in G=GL_2(\Z)/\pm I \, :\, g\omega =\omega\}
\end{equation*}
and $E_\Delta$ was studied. The inverse of $\lambda_\omega$ acts
as
\begin{equation*}
\lambda_\omega^{-1} \bigg( \frac{u+v\sqrt{\Delta}}{2} \bigg)
=\left[ \begin{matrix} \frac{u-Bv}{2} & -Cv \\ Av & \frac{u+Bv}{2}
\end{matrix} \right].
\end{equation*}
It sends $\eps^+_\Delta =Q_\ell \omega+Q_{\ell -1}$ with
$\ell=\eper(\omega)$, to $\Omega^+=\widetilde{M}(\omega)$.
Moreover, one has
\begin{equation*}
\NN\circ \lambda_\omega =\det\qquad \mbox{\rm and} \qquad \tr
\circ \lambda_\omega =\Tr.
\end{equation*}
This implies in particular that $\eps_0 (\omega)$ coincides with
the spectral radius $R( \widetilde{M}(\omega))$ of
$\widetilde{M}(\omega)$, thus
\begin{equation*}
\rho (\omega)=2\log R\big( \widetilde{M}(\omega)\big).
\end{equation*}

Denote by $\RR(\Delta)$ the set of reduced quadratic irrationals
of discriminant $\Delta$ and consider the set
\begin{equation*}
T(N)=\bigg\{ (k,\omega) \, :\, \omega \in \RR=\bigcup_{\Delta >0}
\RR(\Delta),\ \tr \big( \eps^+_\Delta (\omega)^k \big) \leq N
\bigg\}.
\end{equation*}
As seen in Section 2, the latter has cardinality
$\Psi_{\operatorname{ev}}(N)$. As shown in \cite{KOPS}, Proposition 4.3,
the map given by
\begin{equation*}
j(a_1,\dots,a_{2m})=\bigg( \frac{2m}{\eper (\omega)}\, , \omega
\bigg) ,\qquad \mbox{\rm where}\quad \omega
=\overline{[a_1,\dots,a_{2m}]},
\end{equation*}
defines a one-to-one correspondence between $\WW_{\operatorname{ev}}(N)$
and $T(N)$. Denote
\begin{equation*}
r(N)=\sum\limits_{\substack{\omega \in \RR \\ \eps^+ (\omega)<N}}
1 =\pi_0 (2\log N).
\end{equation*}
Then the identification between $\WW_{\operatorname{ev}}(N)$ and $T(N)$
plainly implies as in the proof of \cite{KOPS}, Proposition 4.5,
the inequalities
\begin{equation}\label{5.1}
\sum\limits_{1\leq k\leq 2\log N} r \bigg( \Big(
N-\frac{1}{2}\Big)^{\frac{1}{k}} \bigg)
<\Psi_{\operatorname{ev}}(N)<\sum\limits_{1\leq k\leq 2\log N} r
(N^{\frac{1}{k}}) .
\end{equation}
From the first inequality we infer
\begin{equation}\label{5.2}
r (N)<\Psi_{\operatorname{ev}}(N+1) \ll N^2.
\end{equation}
From \eqref{5.1}, \eqref{5.2}, and Proposition \ref{P1.2} we
derive
\begin{equation}\label{5.3}
\begin{split}
\Psi_{\operatorname{ev}}(N) & <\sum\limits_{1\leq k<2\log N} r
(N^{\frac{1}{k}})<\sum\limits_{1\leq k<2\log N}
\Psi_{\operatorname{ev}}\big( (N+1)^k\big) \\
& =\Psi_{\operatorname{ev}}(N+1)+O \Bigg( \sum\limits_{2\leq k<2\log N}
N^{\frac{2}{k}}\Bigg) \\ & =\Psi_{\operatorname{ev}}(N+1)+O(N\log N) \\
& =\frac{N^2\log 2}{2\zeta (2)} +O_\eps (N^{\frac{7}{4}+\eps}).
\end{split}
\end{equation}
The estimate
\begin{equation*}
r(N)=\frac{N^2\log 2}{2\zeta(2)} +O_\eps (N^{\frac{7}{4}+\eps})
\end{equation*}
follows now immediately from \eqref{5.2} and \eqref{5.3}. This
completes the proof of Proposition \ref{P1.3} by taking
$N=\operatorname{e}^{\frac{X}{2}}$.

{\bf Acknowledgments.} I am grateful to A. Zaharescu for his
comments on a first draft of this paper. I am very thankful to the
referee for criticism that led to a correction of the formula for
$c_2$ and of the estimate for $\Psi(N)$ in Theorem \ref{T1.1},
and for providing a Mathematica notebook which shows some very
interesting features of the quantity $S_N-C_N$.

\end{document}